\documentclass[oneside]{fic-l}
\usepackage{amsmath, amssymb,graphicx,amscd,color,jrsymb}
\newtheorem{thrm}{Theorem}[section]
\newtheorem{prop}[thrm]{Proposition}
\newtheorem{defn}[thrm]{Definition}
\newtheorem*{defn2}{Definition}

\theoremstyle{remark}
\newtheorem*{remark}{Remark:}
\newtheorem*{warning}{Warning:}
\newtheorem*{example}{Example:}
\newtheorem*{remarks}{Remarks:}

\newcommand{\kh}{Kh}
\newcommand{\khr}{Kh_r}
\newcommand{\hfk}{\widehat{HFK}}

\newcommand{\Z}{{\mathbf{Z}}}
\newcommand{\Q}{{\mathbf{Q}}}
\newcommand{\R}{{\mathbf{R}}}
\newcommand{\ts}{{{\thinspace}}}
\newcommand{\rank}{{{\text {rank} \ }}}

\title{Knot Polynomials and Knot Homologies}
\author{Jacob Rasmussen}
\address{ Dept. of Mathematics, Princeton University, \\ Princeton, NJ 08544}
\email{jrasmus@math.princeton.edu}
\thanks{The author was partially supported by an NSF Postdoctoral
  Fellowship.}


\begin{document}

\begin{abstract}
This is an expository paper discussing some parallels between the
Khovanov and knot Floer homologies. We describe the formal similarities
between the theories, and give some examples which illustrate a
somewhat mysterious correspondence between them. 
\end{abstract}

\maketitle

\section{Introduction}

The past few years have seen considerable growth in the theory of what
might be called ``knot homologies.'' Roughly speaking, these
invariants are homological versions of the now--classical knot
polynomials --- the Alexander polynomial, the Jones polynomial, and
their mutual generalization, the HOMFLY polynomial. As a first
approximation, we might say that a knot homology is a bigraded
homology group \( G(K) \) associated to a knot \(K \subset S^3 \). The two
gradings are a ``homological grading'' (the usual sort of grading one
expects on a chain complex) and a ``filtration grading.'' If we take
the filtered Euler characteristic of \(G(K)\), we recover 
 the corresponding knot polynomial.
 
As our understanding of these objects  evolves, it seems likely that the 
definition of a knot homology will evolve with it. As a first
approximation, however, we offer the following:

\begin{defn2}
A knot homology is a theory which assigns to an oriented link \(L \subset
S^3\) together with some auxiliary data \(\mathcal{D}\) 
a filtered chain complex \(\mathcal{C}(L, \mathcal{D})\) satisfying the
following properties:
\begin{enumerate}
\item The filtered Euler characteristic of \(\mathcal{C}(L,
  \mathcal{D}) \) is a ``classical'' polynomial invariant of \(L\). 
\item The filtration on \(\mathcal{C}(L, \mathcal{D})\) gives rise to
  a spectral sequence \(\{E_i, d_i\} \ (i>0) \). For all \(i \geq 2\),
  \(E_i\) does not depend on the choice of auxiliary data \(\mathcal
  D\), and is thus an invariant of the link \(L\). 
\item The homology of the total complex \(\mathcal{C}(L, \mathcal{D})\)
depends only on coarse information about \(L\), such as the number of
components and their linking numbers.
 \end{enumerate}
\end{defn2}
In this formulation, the group \(G(K)\) mentioned above is the \(E_2\)
term of the spectral sequence. 

One thing that makes this definition attractive is the fact that the
known examples arise from rather different areas of
mathematics. The idea that such an object  might exist at all is due to Mikhail
Khovanov. In \cite{Khovanov} he
constructed a bigraded homology theory which I'll call \(Kh(K)\),
whose filtered Euler characteristic is the unnormalized Jones
polynomial of \(K\).  More recently, Khovanov and Rozansky \cite{KHRoz} have
constructed an infinite family of of such knot homologies, one for
each \(n>0\). Their filtered Euler characteristics give
certain specializations of the HOMFLY polynomial. (When \(n=2\) one
recovers \(Kh(K)\).) Our other 
example of a knot homology comes from gauge theory; more precisely,
from the Heegaard Floer homology introduced by Peter Ozsv{\'a}th
and Zolt{\'a}n Szab{\'o} \cite{OS1}. This theory naturally gives rise to a
bigraded homology theory
 known as the knot Floer homology \cite{OS7}, \cite{thesis},
which I'll denote by \(\hfk(K)\). Its filtered Euler characteristic is
the Alexander polynomial of \(K\), which 
 corresponds to the  case \(n=0\) missing from
 Khovanov and Rozansky's construction.

At a first glance, the two types of knot homologies appear to be quite
different. Although they share the formal properties listed above, they are
defined and computed in very different ways, and things which are
easy to see in one
theory may be quite unexpected in the other. For example, it was
obvious from the start that \( \hfk\) is the \(E_2\) term of a
spectral sequence, but the corresponding fact for \(\kh\) was
discovered by Lee in \cite{ESL2}, several years after the appearance of
\cite{Khovanov}. (For the Khovanov-Rozansky theories, this result is due to
Gornik \cite{Gornik}.) On closer inspection, however, a more subtle
correspondence between the two theories begins to appear. This correspondence
has guided much of my own research in this area. In particular,
it led to the discovery of a relation between the Khovanov homology
 and the slice
genus, which was the subject of my talk at McMaster. Rather than
simply rehash this material, which is already covered in \cite{khg}, I
thought I would try to explain where it came from. 

The main goal of this paper, then, is to describe the
above-mentioned correspondence between the Khovanov homology and the
knot Floer homology, and to give some examples to convince
the reader that it is an interesting one. This correspondence does not hold for
 all knots, but it is common enough that the author feels that there
 must be some sort of explanation. Perhaps someone who reads this
 paper will be able to provide one. 
To properly explain the correspondence between the two theories, one
 must summarize a number of  basic facts about them.  As a secondary
 goal, we have tried to make this summary self-contained and 
 accessible to anyone interested in learning about knot homologies.

The rest of the paper is organized as follows. We begin in
section~\ref{Sec:Filt} with a brief review of some facts about
filtered chain complexes. Sections~\ref{Sec:HFK}
and \ref{Sec:Kh} describe the basic properties of the knot
Floer homology and the Khovanov homology, respectively. We do not
attempt to give definitions, but instead focus on the formal
properties of these theories, their relation with classical models for
the Alexander and Jones polynomials, and methods of computation. In
section~\ref{Sec:Correspondence}, we describe the correspondence we
have in mind, and give some reasons for believing that it is
interesting. Finally, in sections~\ref{Sec:Examples} and ~\ref{Sec:Examples2},
 we describe some examples
of knots for which the correspondence is known to hold, as well as
a few cases for which it fails. 

The author would like to thank Matthew Hedden, Peter Kronheimer, 
Ciprian Manolescu, Peter Ozsv{\'a}th, and Zolt{\'a}n Szab{\'o} for many
helpful conversations, and Hans Boden, Ian Hambleton, Andrew Nicas,
and Doug Park for putting together a great conference. 

\section{Preliminaries on Filtered Complexes}
\label{Sec:Filt}
We begin by establishing some notation and conventions related to
filtered chain complexes. Let
 \((\mathcal{C},d)\) be a chain complex freely generated over \(\Z\) by a
 finite set of  generators \(\{x_i\}\). We say that \({\mathcal C}\)
 is a bigraded complex with homogenous generators \(\{x_i\}\) if we
 are given two gradings
\(u: \{x_i\} \to \Z \) and \( f:\{x_i\} \to \Z\)
with the property that if 
\begin{equation*}
 d(x_i) = \sum _j a_{ij}x_j 
\end{equation*}
then \(u(x_j) = u(x_i)-1 \) and \(f(x_j) \leq f(x_i) \) whenever
\(a_{ij} \neq 0 \). We refer to \(u\) as the {\it homological grading}
on \(\mathcal{C}\), and \(f\) as the {\it filtration grading.} 
The filtration grading defines a filtration \(\{\mathcal{F}_n\} \) on
\(\mathcal{C}\), simply by setting
\( \mathcal{F}_n = \text{span} \ts \{x_i \ts | \ts f(x_i) \leq n \}. \)
We refer to \(\mathcal{F}\) as a {\it downward} filtration. (If \(
f(d(x)) \) was always greater than \(f(x) \), the result would be an
{\it upward} filtration.)

\begin{defn}
The filtered Euler characteristic of \(\mathcal{C} \) is
defined to be the sum 
\begin{equation*}
\sum_{i} (-1)^{u(x_i)} f^{f(x_i)}. 
\end{equation*}
It is an element of the Laurent series  ring \( \Z[f] \). 
\end{defn}

The filtration \(\mathcal{F}\) on \( \mathcal{C}\) gives rise to a
spectral sequence, which can be explicitly described as follows. We
 decompose the differential \(d\) in terms of the preferred basis
\(\{x_i\} \), setting
\begin{equation*}
d(x_i) = d_0(x_i) + d_1(x_i) + d_2(x_i)+ \ldots
\end{equation*}
where 
\begin{equation*}
d_n(x_i) = \sum a_{ij}x_j
\end{equation*}
with \(f(x_j) = f(x_i)-n\). Since the number of \(x_i\)'s is
  finite, all but finitely many of the \(d_n\) are \(0\). 
From the identity \(d^2= 0 \), we conclude that \(d_0^2=0 \) as well,
  so \(\mathcal{C}_0 = (\mathcal{C},d_0)\) is a chain complex. Then
  the identity \(d_0d_2 + d_1^2+d_2d_0 = 0 \) implies that 
\( \mathcal{C}_1 = (H(\mathcal{C}_0),d_1)\) is a chain complex.
 Repeating, we obtain a sequence of chain complexes
\( \mathcal{C}_{i+1} = (H(\mathcal{C}_i),d_i)\) which eventually
  converges to \( H(\mathcal{C})\). The complex
  \((\mathcal{C}_i,d_i)\) is generally referred to as the \(E_{i+1}\)
  term of the spectral sequence. 

If we are working with coefficients in a field, it is not difficult to
show ({\it c.f.} section 5.1 of \cite{thesis}) that \(\mathcal{C}_1\)
can be endowed with a differential \(\overline{d}_1\) in such a way
that \( (\mathcal{C}_1, \overline{d}_1)\) is chain homotopy equivalent
to  \( (\mathcal{C},d)\). \( (\mathcal{C}_1, \overline{d}_1)\) is
again a bigraded complex, and the resulting spectral sequence is
isomorphic to our original spectral sequence \(
(\mathcal{C}_i,d_i)\). 

In the definition of the known knot homologies, the following
situation arises. Starting from a link \(L\) plus a choice of some
additional data \(\mathcal{D}\), 
one obtains a bigraded complex \(\mathcal{C}\). A
different choice of data \(\mathcal{D'}\) gives rise to a different
chain complex \( \mathcal{C}'\) together with a filtered chain map
\(\phi: \mathcal{C} \to \mathcal{C}' \). It is easy to see that
\(\phi\) induces maps \( \phi_i : \mathcal{C}_i \to \mathcal{C}'_i
\). One checks directly that the map \( \phi_1\) is an isomorphism; it then
follows from general principles \cite{UsersGuide}
that \(\phi_i\) is an isomorphism for all \(i >0\). The knot homology
is the bigraded complex  \(\mathcal{H}(L)=({\mathcal C}_1, \overline{d}_1)\);
it is well-defined up to filtered isomorphism.

Since \(\mathcal{H}(L)\) will be the primary object of our
attention, we make a few comments specific to it. 
To begin with, observe that \(\mathcal{C}_0\) decomposes as a direct
sum of complexes 
\begin{equation*}
\mathcal{C}_0 = \bigoplus \mathcal{C}_0^j 
\end{equation*}
where \(\mathcal{C}_0^j\) is generated by those \(x_i\) with
\(f(x_i)=j\). As a group, then, we can decompose
\(\mathcal{C}_1 = \bigoplus H(\mathcal{C}_0^j) \). 
When we want to distinguish the summands in \(\mathcal{H}(L)\), we denote
  \( H_i(\mathcal{C}_0^j)\) by  \( \mathcal{H}_i(L,j)\) .

It is often convenient to represent \(\mathcal{H}(L) \) by
its {\it filtered Poincar{\'e} polynomial}:
\begin{defn} The filtered Poincar{\'e} polynomial of \(\mathcal{H}\) is
  given by
{\em
\begin{equation*}
P_{\mathcal{H}(L)}(f,u) = \sum_{i,j} \ts (\rank\mathcal{H}_i(L,j) )
\  u^i f^j.
\end{equation*}
\em}
It is a Laurent polynomial in \(u\) and \(f\). 
\end{defn}
\noindent 
If we substitute \(u=-1\), the filtered Poincar{\'e} polynomial
reduces to the filtered Euler characteristic. When \( \rank
 \mathcal{H}_i(L,j) =1\), we will often use the shorthand \(u^if^j\) to
refer to a generator of this group.

\section{The Knot Floer Homology}
\label{Sec:HFK}

Let \(K\) be a knot in \(S^3\). The knot Floer homology \( \hfk(K)\)
is a bigraded chain complex equipped with a homological
grading \( u\) and a filtration grading \(t\), which is also known as the
 {\it Alexander} grading. Conventionally, the Alexander grading is
 chosen so as to define a downward filtration on \(\hfk(K)\).
   The filtered
Euler characteristic of \( \hfk(K) \) is the Alexander polynomial 
\(\Delta_K(t) \), and the homology of the complex \( \hfk(K)\) is a
single copy of \(\Z\) in homological grading \(0\). 

\begin{figure}
\includegraphics{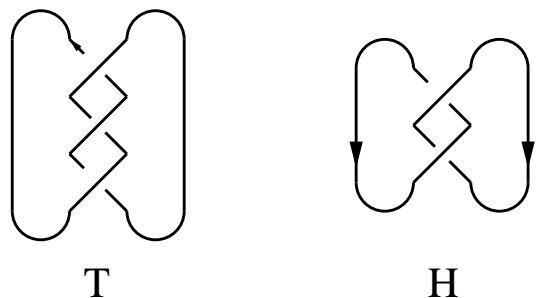}
\caption{\label{Fig:Trefoil} The positive trefoil and Hopf link.}
\end{figure}

\subsection{Examples}
 Let \(U\) be the unknot. The complex 
\(\hfk (U) \) is generated by a single element whose Alexander and
homological gradings are both equal to zero. The differential on this
complex is necessarily trivial. 

Let \(T\) be the positive trefoil knot shown
in Figure~\ref{Fig:Trefoil}. Then  \( \hfk(T)\) is generated by three elements
 \(x_{-1},x_0\), and \(x_{1} \), where
\( t(x_i) = i \) and \( u(x_i) = i-1 \). 
The filtered Poincar{\'e} polynomial is given by 
\begin{equation*}
P_{\hfk(T)}(t,u) = t + t^0u^{-1} + t^{-1}u^{-2}.
\end{equation*}
Substituting \(u=-1\), we recover the filtered Euler characteristic,
 which is the Alexander polynomial of \(T\):
\begin{equation*}
 \Delta_T(t) = t-1 +t^{-1}. 
\end{equation*}
The action of \(d\) is given by \(d(x_{-1}) = d(x_1) = 0\), \(d(x_0) =
x_{-1}\).

\subsection{Symmetry and the \(\delta\)-grading}
It is well known that the Alexander polynomial is symmetric under the
involution which sends \(t \mapsto t^{-1} \). \( \hfk \) is endowed
with an analogous symmetry. To descibe the behavior of the homological
grading under this symmetry, it is convenient to introduce a third
grading on the knot Floer homology.
\begin{defn}
Suppose \( x \) is a homogenous element of \( \hfk(K) \). We define
the \(\delta\)-grading on \( \hfk (K) \) by \( \delta(x) = t(x) -
u(x) \), and denote by \( \hfk(K,i,j)\) the filtered subquotient of
\(\hfk(K) \) with Alexander grading \(i\) and \( \delta\)-grading
\(j\). 
\end{defn}

Since \( t(d(x)) < t(x) \) and \(u(d(x)) = u(d(x))-1\), it follows that
 the \( \delta\)-grading induces a filtration on \( \hfk(K)\). We do
 not really get any new information from this filtration. In fact, it
 is easy to see that the induced spectral sequence is the same as
 the one induced by the Alexander grading. Nonetheless, the 
 \(\delta\)-grading turns out to be
 a convenient and natural thing to consider.  Our first piece of
 evidence for this fact is provided by

\begin{prop}
{\em (\cite{OS7}, \cite{thesis})}
\(\hfk(K,i,j) \cong \hfk(K,-i,j).\)
\end{prop} 
This symmetry is easily seen to hold in the example of the trefoil,
where all the generators have \(\delta\)-grading \(1\).

\subsection{\(\delta\)-thin knots}

Knots for which all generators of \(\hfk\) have the same
\(\delta\)-grading (like the trefoil) have particularly simple knot
Floer homologies.

\begin{defn}
We say that a knot \(K\) is {\em \( \delta\)-thin} with \( \delta =
n\) if all generators of \( \hfk(K) \) have delta-grading \(n\). 
\end{defn}

If \(K\) is \(\delta\)-thin, all generators of \( \hfk(K)\) in a given
Alexander grading have the same homological grading as well. It
follows that the isomorphism class of
 \( \hfk(K)\) is completely determined by the Alexander
polynomial of \(K\) and the value of the \(\delta\)-grading in which
it is supported. 

A large class of \(\delta\)-thin knots is provided by 
\begin{thrm} {\em (\cite{OS8})}
\label{Thm:AlexAlt}
Alternating knots are \(\delta\)-thin with \( \delta =
\sigma(K)/2\). (In  \cite{2bridge} and
\cite{thesis}, such knots were called {\em perfect}.)
\end{thrm}

\begin{warning}
 We use the sign convention of \cite{khg}, namely that 
positive knots  have positive signature. (A {\it positive link} is one
which admits a planar diagram in which all crossings are positive.)
This is the opposite of the convention
used in \cite{OS8}.
\end{warning}

Many small nonalternating knots are \( \delta\)-thin as well. Among
the 53 nonalternating knots with 10 or fewer crossings, at least 39
are known to be \(\delta\)-thin. The simplest example of a knot which
is not \(\delta\)-thin is the \((3,4)\) torus knot, whose Poincar{\'e}
polynomial is given by
\begin{equation*}
P_{\hfk(T_{3,4})} (t,u) = t^3 + t^2u^{-1} + u^{-2} + t^{-2}u^{-5} +
t^{-3}u^{-6}
\end{equation*}

\subsection{Relation with the Knot Genus}

 \(\hfk (K) \) carries a lot
of geometric information about the knot \(K\) and the manifolds
obtained by surgery on it. In particular, it provides lower bounds on the
genus of embedded surfaces bounding \(K\), both in \(3\) and \(4\)
dimensions. In three dimensions, this bound  is actually
sharp:
\begin{thrm}
\cite{OSGen}
Let \(g(K)\) denote the Seifert genus of \(K\). Then for all \( i >
g(K) \), \(\hfk(K,i) = 0\), while \( \hfk(K,g(K)) \neq 0 \). 
\end{thrm}
This generalizes the well-known fact that \(g(K)\) is greater
than or equal to the degree of \(\Delta_K(t)\). 

\subsection{The \(\tau\)-invariant}

The Alexander filtration on \(\hfk (K)\) gives rise to a spectral
sequence, all of whose terms are invariants of \(K\). This spectral
sequence converges to the homology of the total complex, which is \(\Z\). 

\begin{defn}
{\em (\cite{OS10}, \cite{thesis})}
Let \(\tau(K) \) be the Alexander grading of the surviving copy of
\(\Z\) in the spectral sequence for \(\hfk(K)\). 
\end{defn}

Since the spectral sequence is an invariant of \(K\), \(\tau\) is
clearly an invariant as well. Below, we summarize some interesting
properties of \( \tau\):

\begin{prop}
\label{Prop:TauProp}
The invariant \( \tau(K) \) satisfies the following:
\begin{enumerate}
\item (Additivity) \(\tau(K_1\#K_2) = \tau(K_1) + \tau (K_2) \). If
  \(\overline{K} \) is the mirror image of \(K\), then
  \(\tau(\overline{K}) = -\tau(K) \). 
\item (Adjunction) \( |\tau(K)| \leq g_*(K)\), where \(g_*(K) \)
  denotes the slice genus of \(K\). 
\item If \(K\) is an alternating knot, then \(\tau(K) = \sigma (K) /2
  \). 
\item If \(K\) is a positive knot, then 
\( \tau(K) = g_*(K) = g(K) \). 
\end{enumerate}
\end{prop}

\begin{remarks} Properties (1)--(3) are due to Ozsv{\'a}th
and Szab{\'o}, and may be found in \cite{OS10}.
(Property (3) is a corollary of Proposition~\ref{Thm:AlexAlt}.)
 They also proved
property (4) for the special case of torus knots \cite{OSLens}. The
general case follows from this special one, together with  work
of Livingston \cite{Livingston} and Rudolph \cite{Rudolph}. 

Property (2) is an application of the adjunction inequality in
Ozsv{\'a}th-Szab{\'o} theory. This inequality is familiar from
classical gauge theory, and was first applied in this context by
Kronheimer and Mrowka in their proof of the Milnor conjecture \cite{KMMilnor}.
\end{remarks}

\subsection{Links and the skein exact sequence}

Let \(L \subset S^3\) be an oriented \(n\)-component link. In
\cite{OS8}, Ozsv{\'a}th and Szab{\'o} show how \(L\) can naturally be
thought of as a knot in \(\#^{n-1}(S^1\times S^2) \). This
construction gives rise to a knot Floer homology group \(
\hfk(L) \), which is again a filtered complex. Its 
 filtered Euler characteristic is given by 
\begin{equation*}
P_{\hfk(L)}(t,-1) = (t^{1/2}-t^{-1/2})^{n-1} \Delta_L(t),
\end{equation*}
 and its total homology has rank \(2^{n-1}\). The 
Poincar{\'e} polynomial of the total homology is given by
\begin{equation*}
P(u) = (u^{1/2} + u^{-1/2})^{n-1}.
\end{equation*}
(when \(n\) is odd, the homological grading on \( \hfk(L) \) is naturally
an element of \( \Z + \frac{1}{2}\) rather than of \(\Z\).)


\begin{prop}
\(\hfk(L) \) has the following elementary properties:
\begin{enumerate}
\item \( \hfk(L^o) \cong \hfk(L) \), where \( L^o\) denotes \(L\) with
  the orientations of all components reversed.
\item \( \hfk(\overline{L}) \cong \hfk(L)^* \), where \(
  \overline{L}\) is the mirror image of \(L\), and \(^*\) denotes the
  operation of taking the dual complex.
\item \( \hfk(L_1\#L_2) \cong  \hfk(L_1) \otimes \hfk(L_2) \), where
  \( L_1 \# L_2\) is the link obtained by taking the oriented
  connected sum of any component of  \(L_1\) with any component of \(L_2\).
\item \(\hfk(L_1 \coprod L_2) \cong \hfk(L_1)\otimes \hfk(L_2) \otimes
  X\), where \(X\) is the rank two complex with Poincar{\'e}
  polynomial \(P_X(t,u) = u^{-1/2} + u^{1/2} \) and trivial differential. 
\end{enumerate}
\end{prop}

In addition, \(\hfk\)  satisfies a {\it skein exact
  sequence}, which is a generalization of the skein relation for the
  Alexander polynomial:

\begin{prop}
\label{Prop:SkeinAlex}
{\em (\cite{OS8})} There are  long exact sequences
\begin{equation*}
\begin{CD}
@>>> \hfk_{*-1/2}(\undercrossing) @>>> \hfk_*(\asmoothing) @>>>
\hfk_{*-1}(\overcrossing) @>>>
\end{CD}
\end{equation*}
(when the middle term has more components than the other two terms)
and
\begin{equation*}
\begin{CD}
@>>> \hfk_{*-1/2}(\undercrossing)  @>>> \hfk_*(\asmoothing ) \otimes A @>>>
\hfk_{*-1}(\overcrossing) @>>>
\end{CD}
\end{equation*}
(when the middle term has fewer components.) Here \(A\) is the complex
with filtered Poincar{\'e} polynomial \(P_A =
[(tu)^{-1/2}+(tu)^{1/2}]^2 \) and trivial differential.
 All the maps in these sequences respect the Alexander filtration. 
\end{prop}

A generalization of Theorem~\ref{Thm:AlexAlt} holds as well: {\it
  nonsplit} alternating links are \(\delta\)-thin \cite{OS8}. 

\begin{example} Let \(H\) denote the positive Hopf link of
Figure~\ref{Fig:Trefoil}. \(\hfk (H) \) is free of rank 4, and its
filtered Poincar{\'e} polynomial is given by 
\begin{equation*}
P_{\hfk (H)} = tu^{1/2} + 2 u^{-1/2} + t^{-1}u^{-3/2}.
\end{equation*}
In the skein exact sequence
\begin{equation*}
\begin{CD}
@>>> \hfk_*(U ,i) @>>> \hfk_{*-1/2}(H ,i) @>>>
\hfk_{*-1/2}(T ,i) @>{f}>>
\end{CD}
\end{equation*}
the map \(f\) is the \(0\) map.
\end{example}

\subsection{Methods of Computation}

Recall that
\( \hfk (K) \) arises as the second term in the spectral sequence
of a certain bigraded chain complex, which we will call
\(\widehat{CFK}(K) \). The extra data needed to define \(
\widehat{CFK}(K)\) is a Heegaard splitting of the complement of \(K\)
together with a preferred meridian for this splitting. Given this
data, the generators of \(\widehat{CFK}(K) \)  can be computed by 
a process which is more or less the same as computing the Alexander
polynomial via Fox calculus \cite{thesis}. 
In contrast, the differentials in the complex
\(\widehat{CFK}(K) \) are determined by  counting the number of
elements in 
 certain zero-dimensional moduli spaces of pseudoholomorphic
disks in a symplectic manifold. Although it is sometimes possible to
determine the number of points in these moduli spaces, it is in general a
difficult problem. As a result, there is currently no known algorithm
for computing the knot Floer homology of a given knot. 

In most cases, successful computation of \( \hfk\) depends on finding
a nice Heegaard splitting for the knot complement. Two
particularly nice classes of splittings were described by Ozsv{\'a}th
and Szab{\'o} in \cite{OS7} and \cite{OS8}. The first type of
splitting is applicable to a specific class of knots --- those which
can be represented by a doubly pointed Heegaard diagram of genus
\(1\). In \cite{GMM}, Goda, Matsuda, and Morifuji observed that these
knots are precisely those which admit \((1,1)\) bridge
decompositions. (See section~\ref{Subsec:(1,1)} for more details.) 
For such knots, the methods of \cite{OS7} provide a
completely algorithmic way of computing the knot Floer homology. 

The method of \cite{OS8} is based on the Kauffman state model for the
Alexander polynomial \cite{Kauffman2}. It is potentially
applicable to any knot, but is most effective for alternating knots
--- it is used to prove Theorem~\ref{Thm:AlexAlt} --- and for knots with
relatively small crossing number. It has been used by Ozsv{\'a}th,
Szab{\'o} \cite{OSMut}, and Eftekhary \cite{Eaman1} to compute the knot
Floer homology of three-strand pretzel knots. As a rule of thumb, it
tends to be effective at computing \( \hfk(K,i)\) when \(i\) is close
to \(g(K)\), but rather less so when \(i\) is close to \(0\). 
Other special Heegaard splittings have been used by Eftekhary
\cite{Eaman2} and Hedden \cite{Hedden} to make some computations for
Whitehead doubles and cabled knots. 

We close this section by mentioning two indirect
computational techniques, which do not rely on a choice of
Heegaard splitting. The first is to use the skein exact sequence of 
Proposition~\ref{Prop:SkeinAlex}. This can be an effective method for
proving a knot is \(\delta\)-thin, especially when the knot is only
mildly nonalternating or has a small number of crossings. The second
method applies if \(K\) has a lens space, or, more generally an
\(L\)-space surgery. In this case, \(\hfk(K)\) is completely
determined by the Alexander polynomial of \(K\) \cite{OSLens}. 

\section{The Khovanov Homology}
\label{Sec:Kh}

Let \(L \subset S^3\)  be an oriented \(n\)--component link. The
Khovanov homology \( \kh (L) \) is a bigraded chain complex equipped
with  a homological grading \(u \) and a filtration grading \(q\),
which is also known as {\it
  Jones} grading. As a group, \( \kh (L) \) was defined by
Khovanov in \cite{Khovanov}; the chain complex structure was described
by  Lee in \cite{ESL2}. Conventionally, \(\kh\) is defined to be a
{\it cohomology} theory with an upward filtration. 
The filtered Euler characteristic of \( \kh (L) \) is given by
\begin{equation*}
P_{\kh(L)} (q,-1) = (q+q^{-1})V_L(q^2),
\end{equation*}
where \(V_L(t)\) is the Jones polynomial of \(L\).
 
The homology of the
complex \( \kh(L) \) has rational rank \(2^n\) and has no
\(p\)--torsion for \(p\neq 2\), but its \(2\)--torsion can be rather
complicated.
If \(n=1\) (so \(L\) is actually a knot) then both rational 
generators have homological grading \(0\). More generally, for
\(n>1\), the homological gradings of the generators are determined by
the pairwise linking numbers of the components of \(L\). (See
\cite{ESL2} for details.) 

Some elementary properties of the Khovanov homology are stated
below:
\begin{prop} \(\kh(L) \) satisfies 
\begin{enumerate}
\item \( \kh(L^o) \cong \kh(L) \).
\item \( \kh(\overline{L}) \cong \kh(L)^* \).
\item \( \kh(L_1 \coprod L_2) \cong  \kh(L_1) \otimes \kh(L_2) \).
\end{enumerate}
\end{prop}

\subsection{Examples} Let \(U\) be the unknot. \( \kh(U) \) is
generated by two elements \(x_{\pm}\). Both generators have
homological grading 0, and their \(q\)-grading is given by
\(q(x_{\pm}) = \pm 1 \). The graded Poincar{\'e} polynomial is
\begin{equation*}
P_{\kh(U)} (q,u) = q^{-1} + q.
\end{equation*}
The differential on \(\kh(U)\) is necessarily trivial. 

Let \(H\) be the Hopf link of Figure~\ref{Fig:Trefoil}. Then
\( \kh(H) \) has rank \(4\), and its graded Poincar{\'e} polynomial
is 
\begin{equation*}
P_{\kh(H)} = 1 + q^2 + q^4u^2 + q^6u^2.
\end{equation*}
As a complex \(\kh(H) \) is trivial, so its homology has rank
\(4\). 

Let \(T\) be the trefoil of
Figure~\ref{Fig:Trefoil}. Then \(\kh (T) \) has rational rank \(4\),
and its graded Poincar{\'e} polynomial is 
\begin{equation*}
P_{\kh(T)}(q,u) = q + q^3 + q^5u^2 + q^9 u^3.
\end{equation*}
If we use \(\Z/2\) coefficients, however, the rank is \(6\):
\begin{equation*}
P_{\kh(T;\Z/2)}(q,u) = q + q^3 + q^5u^2 + q^7u^2 + q^7 u^3+  q^9 u^3.
\end{equation*}
We have 
\begin{align*}
P_{\kh(T)}(q,-1) & = q+ q^3 + q^5 -q^9 \\
& = (q^{-1}+ q)(q^2 + q^6-q^8) \\
& = (q^{-1}+ q)V_T(q^2). 
\end{align*}
There is a single nonzero differential in the complex \(\kh (T) \),
which takes \(q^5u^2\) to \(q^9u^3\). The total homology thus has rank \(2\),
with both generators having homological grading zero. 

\subsection{The skein exact sequence}


The Khovanov homology is constructed using the Kauffman state model
for the Jones polynomial. As such, it is naturally endowed with a
skein exact sequence based on Kauffman's {\it unoriented} skein
relation for the Jones polynomial. 

\begin{prop}
\label{Prop:KhSkein}
There are long exact sequences
\begin{equation*}
\begin{CD}
@>{\cdot u}>> q^{2+3\epsilon} u^{1+\epsilon} \kh(\hsmoothing) @>>> \kh(\overcrossing)
@>>> q \kh (\asmoothing) @>{\cdot u}>>
\end{CD}
\end{equation*}
and 
\begin{equation*}
\begin{CD}
@>{\cdot u}>> q^{-1} \kh(\asmoothing) @>>> \kh(\undercrossing)
@>>>q^{1+3\epsilon}u^\epsilon
 \kh (\hsmoothing)q^{-1}  @>{\cdot u}>>
\end{CD}
\end{equation*}
where \(\epsilon \) is the
difference between the number of negative crossings in the unoriented
resolution \(\hsmoothing \) and the number of such crossings 
in the original diagram. 
\end{prop}
\noindent Here, notation such as 
\( q \kh (\asmoothing)\) should be understood to
indicate the complex \(\kh (\asmoothing)\)
shifted in such a way as to multiply its Poincar{\'e} polynomial
by \(q\). The arrow marked with \(\cdot u \) is the boundary map in
the long exact sequence; it raises the homological grading by \(1\).

\begin{example} Let \(H\) be the positive Hopf link. Then both
resolutions of a crossing yield the unknot, and the first exact
sequence becomes
\begin{equation*}
\begin{CD}
@>>> q^5u^2 \kh(U) @>>> \kh(H) @>>> q\kh(U) @>>>
\end{CD}
\end{equation*}
which splits to give a short exact sequence
\begin{equation*}
\begin{CD}
0 @>>> q^5u^2(q+q^{-1}) @>>> \kh(H) @>>> q (q+q^{-1}) @>>> 0. 
\end{CD}
\end{equation*}
\end{example}

Just as the unoriented skein relation for the Jones polynomial can be
used to show that it satisfies the oriented skein relation

\begin{equation*}
q^{-2}V_L(\overcrossing)-q^2V_L(\undercrossing) =
(q-q^{-1})V_L(\asmoothing)\
\end{equation*}
 the skein exact sequence above can
be used to show that \(\kh\) satisfies an oriented skein exact
sequence analogous to that of Proposition \ref{Prop:SkeinAlex}. 

\subsection{The reduced Khovanov homology}

In section 3 of \cite{Khovanov2}, Khovanov describes a slight variant of his
construction which results in  a related bigraded homology
theory known as the {\it reduced Khovanov homology.} This group is an
invariant of a link \(L\) together with a particular  marked
 component \(L_i\) of \(L\). We denote it by  \( \kh_r(L,L_i) \), or
  \( \khr(K)\) if \(K\) is a knot. Like
 \(\kh \), \( \kh_r\) is endowed with a homological grading \( u\)
 and a Jones grading \(q\). Its graded  Euler characteristic
 is given by the Jones polynomial:

\begin{equation*}
P_{\khr (L,L_i)} (q, -1) = V_L(q^2). 
\end{equation*}

Recent work of Bar-Natan \cite{DBN2} and Turner \cite{Turner} implies
that \(\kh_r\) may be endowed with a differential analogous to Lee's,
but without the problems with \(2\)-torsion. The total homology of
the complex \(\kh_r\) is \(\Z^{n-1}\), where \(n\) is the number of
components of \(L\). When \(L\) has more than one component, \(\kh_r\)
suffers from the disadvantage that it is not really a link invariant:
it depends on the choice of marked component. For knots, however,
\(\khr \) seems to be 
an interesting and natural invariant in its own right. 

The reduced Khovanov homology satisfies skein exact sequences
analogous to the ones described above for \( \kh\). 
The two theories are  related by 
\begin{prop}
There is a long exact sequence 
\begin{equation*}
\begin{CD}
q^{-1} \kh_{r  *} (L,L_i) @>>> \kh _*(L) @>>> 
q \kh_{r*} (L,L_i) @>{\partial_r}>> q^{-1} \kh_{r  *+1} (L,L_i)
\end{CD}
\end{equation*}
\end{prop}

\begin{example}
 Let \(T\) be the positive trefoil knot. Then 
\(\khr(T) \cong \Z^3 \); its graded Poincar{\'e} polynomial is 
\begin{equation*}
P_{\khr (T)}(q,u) = q^2 + q^6t^2 + q^8t^3.
\end{equation*}
The boundary map \( \partial_r\) takes \(q \khr(T)\), which has
Poincar{\'e} polynomial  
\begin{equation*}
q^3 + q^7t^2 + q^9t^3 
\end{equation*}
 to 
\(q ^{-1} \khr(T) \), which has Poincar{\'e} polynomial
\begin{equation*}
q + q^5t^2 + q^7t^3.
\end{equation*}
 The only possible nontrival component of \(\partial_r\) is the one
 which takes \(q^7t^2\) to   \(q^7t^3\). It is not difficult to
 see that this component is multiplication by \(2\). 
\end{example}

\begin{remark}
 In fact, it follows from work of Ozsv{\'a}th
and Szab{\'o} \cite{OSKh} or Shumakovitch \cite{Shumak2} that the map
\( \partial_r\) is always congruent to \(0\) \(\text{mod} \ 2\). 
\end{remark}

\subsection{The \(\delta\) grading and alternating links}

Just as in the case of the knot Floer homology, it turns out that
there is a third interesting grading on \(\kh\) and \(\khr\).

\begin{defn}
Suppose \(x\) is a homogenous element of \( \kh(L)\). We define the \(
\delta\)-grading on \(\kh(L)\) by \( \delta (x) = q(x) - 2 u(x)
\), and denote by \( \kh ( L, i,j) \) the filtered subquotient of \(
\kh (L) \) with Jones grading \(i\) and \( \delta\)-grading \(j\). 
\end{defn}

If \(K\) is a knot, the  \( \delta \)-grading on \( \khr (K) \) is
defined similarly. 
\(\delta\)-gradings of elements of \( \kh(K)\) are
always odd, while \(\delta\)-gradings of elements of \( \khr (K) \)
are always even. 

\begin{defn}
A knot \(K\) is \(\kh\)-thin with \( \delta = n\) if \( \khr(K) \) is
free over \(\Z\) and all its generators have \(\delta\)-grading
\(n\).
\end{defn}

It follows that  \(\khr\) of a \(\kh\)--thin knot is determined by its
Jones polynomial and the value of \(\delta\) that \(\khr\) is
 supported in. In fact, using the chain complex structure on \(\kh
 (K)\), Lee has shown that the rational \( \kh \) of a \(\kh\)-thin knot is
 determined by this information as well. 

As evidence that the preceding two definitions are interesting, we
have the following theorem, which is essentially due to Lee, although
she phrased it in terms of \(\kh\) rather than \( \khr\). 

\begin{thrm} {\em (\cite{Lee})}
\label{Thm:JonesAlt}
If \(L\) is a nonsplit alternating link, then \(L\) is \(\kh\)-thin
with \(\delta = \sigma (L) \). 
\end{thrm}

\subsection{The \(s\)-invariant}
Let \(K\) be a knot.
If we use rational coefficients, the spectral sequence induced on \(
 \kh(K) \) by the Jones filtration converges to \( \Q^2\). 
 In analogy with  the \(\tau\) invariant, we can define invariants of
 \(K\) by looking at the \(q\)-gradings of the surviving terms in the
 spectral sequence. At first glance, it appears that we get two such
 invariants, since there are two surviving generators in the spectral
 sequence. In reality, it can be shown that the \(q\)-gradings of these 
 generators always differ by \(2\), so there is really only one
 invariant:

\begin{defn}
If \(K\) is a knot in \(S^3\), we let \(s(K)\) be the average of the
\(q\) gradings of the two surviving rational generators in the
spectral sequence for \(\kh(K) \). 
\end{defn}
\noindent Since the two \(q\) gradings are odd integers, \(s(K)\) is even.
The invariant \(s\) can be shown to have the following properties,
which are exact analogs of the properties of \(\tau\) described in
Proposition~\ref{Prop:TauProp}.

\begin{prop} {\em (\cite{khg})}
\label{Prop:sProp}
The invariant \( s(K) \) satisfies the following:
\begin{enumerate}
\item (Additivity) \(s(K_1\#K_2) = s(K_1) + s (K_2) \). If
  \(\overline{K} \) is the mirror image of \(K\), then
  \(s(\overline{K}) = -s(K) \). 
\item \( |s(K)| \leq 2 g_*(K)\), where \(g_*(K) \)
  denotes the slice genus of \(K\). 
\item If \(K\) is an alternating knot, then \(s(K) = \sigma (K)
  \). 
\item If \(K\) is a positive knot then 
\( s(K) = 2 g_*(K) = 2 g(K) \). 
\end{enumerate}
\end{prop}

\subsection{Methods of computation}

 In \cite{Khovanov}, \(\kh(L) \) is defined as the homology of a
 finite dimensional chain complex \(CKh(L)\), which is defined using
 a planar diagram of \(L\). The generators of \( CKh (L) \)
 correspond to states in the Kauffman state model for the Jones
 polynomial. The differentials are completely explicit as well, so 
\( \kh (L) \) is by definition algorithmically computable.
The size of \(CKh (L) \)
 grows exponentially with the number of crossings in the planar
 diagram of \(L\),  so it is only in the simplest cases that
the homology can be computed directly by hand. 
On the other hand, the problem is well-suited to computer
 computation. The first  program for this purpose was written by Bar-Natan
 \cite{DBN}; it could be used to compute \(\kh(L)\) for links of up
 to \(12\) or \(13\) crossings. Based on his calculations, Bar-Natan
 formulated some influential conjectures, which formed the
 basis for much of the early work on Khovanov homology.
More recently, Shumakovitch has written a substantially faster program known as
{\it KhoHo} \cite{Shumak} which can effectively compute \(\kh \) and
 \(\khr \) for links with as many as \(20\) crossings. Despite this
 fact, many basic computational questions remain unanswered. For
 example,  the Khovanov homology of the \((p,q)\) torus knot is
 still unknown. 


\section{The FK Correspondence}
\label{Sec:Correspondence}

We are now in a position to describe the correspondence alluded to in
the introduction. First, we need one more bit of notation. We denote
by \( \hfk(K,*,j) \) the group generated by all generators of \(
\hfk (K)\) which have \(\delta\)-grading equal to \(j\). In other
words, we have 
\begin{equation*}
\hfk(K,*,j) = \bigoplus _{i \in \Z} \hfk (K,i,j).
\end{equation*}
The notation \( \khr (K,*,j) \) should be understood similarly.

\begin{defn}
We say that a knot \(K\) has {\em property FK} (for Floer-Khovanov) if
it satisfies the following two conditions:
\begin{enumerate}
\item For each value of \(j\), we have
\begin{equation*}
\rank \hfk (K,*,j) = \rank \khr (K,*,2j). 
\end{equation*}
\item \( s(K) = 2 \tau (K) \). 
\end{enumerate}
\end{defn}

The definition is motivated by the following

\begin{prop}
Alternating knots have property FK.
\end{prop}
\begin{proof}
This follows trivially from Theorems~\ref{Thm:AlexAlt} and
\ref{Thm:JonesAlt}, since for an alternating knot
\begin{equation*}
\rank \hfk (K)  = |\Delta_K(-1)| = |\det K| = | V_K(-1)| = \rank \khr (K)
\end{equation*}
\end{proof}

What is perhaps more surprising is that a great many non-alternating
knots have property FK as well. In fact, I spent about
six months under the impression that property FK might hold for
all knots before discovering that part (1) of the property
 fails for the \((4,5)\) torus
knot. Part (2) still holds in this case, and to the best of my
knowledge, there still no examples known for which it fails. 

 In section~\ref{Sec:Examples}, we describe several examples of
knots known to have property FK. Many have Floer homologies and
Khovanov homologies which seem quite nontrivial. Although it now seems
likely that property FK fails for all knots which are sufficiently
complicated in some sense, the level of complication needed seems
rather high. When one considers how different the definitions of the
two theories seem, the correspondence seems to demand an explanation. 

\subsection{Possible explanations}
In this section, we describe a number of arguments which might be
advanced to explain the correspondence described above. Although they
are all at least vaguely plausible, none of them seem truly
satisfactory.
\begin{itemize}
\item \textbf{Skein Theory:} The two theories share some similar basic
  properties. They agree for alternating knots and links, both satisfy
  skein exact sequences (If one wants, the unoriented skein sequence
  of Proposition~\ref{Prop:KhSkein} can be used to prove an oriented
  skein sequence similar to that of Proposition~\ref{Prop:SkeinAlex}.)
  and are both constrained by the requirement that they be a complex
  with simple homology. Perhaps these requirements are enough to
  force Property FK to hold for a large number of knots. The arguments
  against this idea are twofold. First, the requirements described
  above are fairly weak in practice. Skein theory can be used to prove
  some things, but there are many knots for which it seems
  ineffective. Second, the Floer homology of the branched double cover
  satisfies at 
  least the first two of the properties described above, but 
  quickly begins to differ from them once you leave the realm of
  alternating knots.

\item\textbf{A Master Theory:} It is tempting to imagine that the two
  theories can be subsumed as special cases of a single construction,
  like that of Khovanov and Rozansky \cite{KHRoz}. The similarity
  might become evident in this more complicated theory. The major
  objection to this idea is that the Khovanov-Rozansky theories do not
  share the simple behavior that \(\khr \) and \(\hfk\) exhibit for
  alternating knots.
  (This follows from the fact that the HOMFLY
  polynomial of an alternating knot need not be alternating.)

\item\textbf{A Spectral Sequence:} A third possibility is that the two
  theories are related by a spectral sequence. This is particularly
  attractive in light of the work in \cite{OSKh}, which showed that
  there is a spectral sequence starting at the reduced
Khovanov homology and converging to 
the Floer homology of the branched double cover. It is also supported
  by the fact that in all the known examples of knots which do not
  have property FK, the rank of the \(\khr\) is greater than that of
  \( \hfk\). This is perhaps the most attractive possibility of the
  three, but it is not clear where such a spectral sequence might come
  from. 
\end{itemize}

\subsection{Applications}

Whatever its origins, the correspondence between the Floer homology
and the Khovanov homology has proved to be a useful guide to the study
of both. The two have complementary strengths and weaknesses. On the
one hand, the Khovanov homology is very simple and easy to compute
with, but we have relatively little geometric intuition into its
behavior. On the other, the Floer homology can be difficult to
compute, but comes with twenty years worth of geometric intuition
developed by gauge theory. 

Even though we have no direct evidence of a relation between the two
theories, the FK correspondence can be a useful guide, suggesting
that we try to prove  analogs of statements which are known in one
theory directly in the other.
For example, the possibility
that alternating knots might be \(\delta\)-thin was first suggested by
Lee's proof of the analogous result for the Khovanov
homology.

Conversely, the fact that \(\tau\) was known to be a lower
bound for the slice genus suggested that one should try to prove the
same thing for its counterpart \(s\), in the Khovanov theory
\cite{khg}. As a corollary, one obtains 
topological proofs of some results which were previously only known
via gauge theory. The classical Milnor conjecture, which states that
the slice genus of a torus knot is equal to its Seifert genus, is
one such case. (Since torus knots are positive, the result is a
consequence of the third property in Proposition~\ref{Prop:sProp}.)

While we were at McMaster, Bob Gompf kindly pointed out another such
application. Namely, \(s\) can also be used to give a gauge-theory free
proof of the existance of an exotic \(\R^4\). Indeed, Gompf has shown
that to construct such a
manifold, it suffices to exhibit a knot \(K\) which is smoothly but
not topologically slice.
 (See \cite{GompfStip} p. 522 for a proof.) By a theorem of Freedman,
any knot with Alexander polynomial \(1\) is topologically
slice \cite{Freedman},
 so we need only find a knot \(K\) with \(\Delta_K(t) = 1\) and \(
s(K) \neq 0 \). It is not difficult to produce such a knot --- for
example, the \((-3,5,7)\) pretzel knot will do. 
The Khovanov homology of this knot can be calculated, 
either by {\it KhoHo} or using the skein exact sequence, and from
there it can easily be determined that \(s=-1\).

\section{Knots with Property FK}
\label{Sec:Examples}
In this section, we describe some examples of knots which
can be seen to have property FK. 
Our main goal is to convince the reader that this property is
both interesting and common. To this end, we have tried to  describe
a variety of knots for which it
holds, including some for which \(\khr \) and \( \hfk \) seem quite
complicated.

\subsection{Knots with few crossings}
In this category we include all knots with 10 or fewer crossings,
numbered as in Rolfsen \cite{Rolfsen}. It seems very likely that all
of these knots have property FK. More precisely, \(\hfk\) of the
non-alternating \(8\) and \(9\) crossing knots was computed by Ozsv{\'a}th
and Szab{\'o} in \cite{OS8}. All of these knots are \(\delta\)-thin
except for \(8_{19}\) (the \((3,4)\) torus knot) and \(9_{42}\). 
For the
10 crossing knots, Goda, Matsuda, and Morifuji computed \(\hfk \) for
\(10_{124}\)--\(10_{139}\), \(10_{145}\) and \(10_{161}\), using the
fact that these are \((1,1)\) knots \cite{GMM}. Of the remainder,
skein theory can be used to show that \(10_{140}\),\(10_{143}\),
\(10_{144}\), \(10_{146}\)--\(10_{149}\),\(10_{151}\),\(10_{155}\),
\(10_{158}\), and \(10_{163}\)--\(10_{166}\) are \(\delta\)-thin with
\(\delta = \sigma /2 \). Finally, the methods of \cite{OSMut} can
be used to show that \(10_{141}\),
\(10_{156}\),\(10_{157}\),\(10_{159}\), and \(10_{160}\) are \(\delta\)-thin
with \( \delta = \sigma/2\), and to compute \(\hfk\) of \(10_{152}\),
\(10_{153}\) and \(10_{154}\), which are not \( \delta\)-thin. 
This leaves two knots --- \(10_{141}\) and
\(10_{150}\) --- for which the author was unable to determine \(
\hfk\). 

The Khovanov homology of all these knots can be computed using either
Bar-Natan's program or {\it KhoHo}, and in all cases it is easy to see
that property FK holds. Perhaps the most interesting example is
provided by the knot \(10_{145}\), for which the ranks of both \( \khr\)
and \(\hfk\) are equal to \(13\). In both cases, this knot exhibits a lot of
``hidden'' homology --- the sum of the absolute values of
coefficients of the Alexander polynomial is only 7, and the
corresponding sum for the Jones polynomial is only 5. Of course, it
is easy to compute the Khovanov homology for \(10_{141}\) and
\(10_{150}\) as well --- they are both \(Kh\)-thin with \(\delta =
\sigma/2\). It would be very surprising if they were not \(
\delta\)-thin too.


\subsection{\((1,1)\) knots}
\label{Subsec:(1,1)} 
A knot \(K \subset S^3\) is said to be 
a \((1,1)\) knot if there is a genus 
\(1\) Heegaard splitting \(S^3 = H_1 \cup_{T^2} H_2 \) of \(S^3\) 
 with the property that \(K \cap H_i\) is a single trivially embedded arc. 
The knot Floer homology of these knots is algorithmically computable 
\cite{OS8}, \cite{GMM}, so they offer us a good opportunity to
compare \(\hfk\) and \(\khr\) on a large set of ``complicated'' knots.
By looking at examples of this type, we were able to turn up a small
sample of knots which do not have property FK, as well as a rather
larger number of knots that do.  
We briefly sketch the method of computation here.

\begin{figure}
\includegraphics{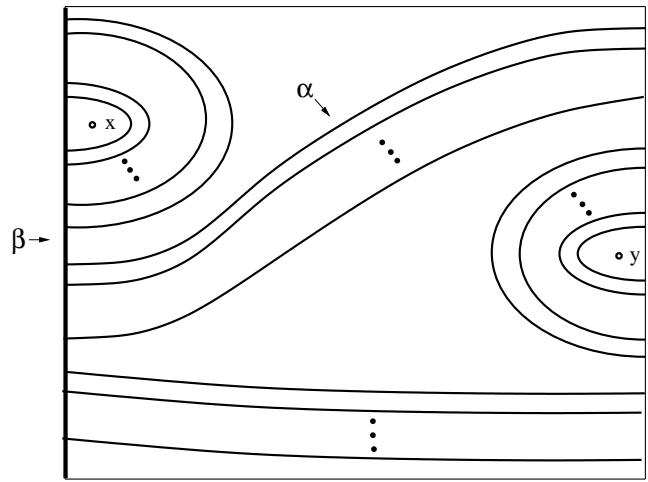}
\caption{\label{Fig:(1,1)Knot} A doubly pointed genus one Heegaard
  diagram for the knot \(K(p,q,r,s)\). The ``rainbow'' on either side
  contains \(q\) strands of \(\alpha\), the middle band contains
  \(r\), and the lower band contains \(p-2q-r\).}
\end{figure}

As described in \cite{GMM}, a \((1,1)\) decomposition of a knot
determines and is determined by a doubly pointed diagram of \(S^3\).
 Such a diagram is composed of 
 a pair of curves \( \alpha\) and \(\beta \) on the
doubly punctured torus \(T^2 - x - y\), with the property that the
algebraic intersection number \(\alpha \cdot \beta = \pm 1 \). After
an isotopy, we may assume that \( \alpha \) and \( \beta \) are in the
form shown in Figure~\ref{Fig:(1,1)Knot}. Such a diagram is 
specified by  four non-negative integers \(p,q,r,\) and \(s\). Here
\(p\) is the total number of intersection points of \( \alpha \) with
\( \beta \), \(q\) is the number of strands in each ``rainbow,'' \(r\)
is the number of strands running from below the left-hand rainbow to
above the right-hand one, and \(s\) is the ``twist parameter'': if we
label the intersection points on either side of the diagram starting
from the top, then the \(i\)-th point on the right-hand side is
identified with the \((i+s)\)-th point on the left-hand side.

Conversely, suppose we are 
given \( p,q, r, s \geq 0 \) satisfying \(2q+r \leq p\)
and \(s < p \), and  with the property that the resulting curve \( \alpha \)
 has intersection number \( \pm 1 \) with \( \beta \). Then we can
recover the corresponding knot in \(S^3\) as follows. First, let \(S^3
= H_1 \cup H_2 \) be the standard genus one Heegaard splitting of
\(S^3\), and let \(\gamma_i\) be a curve on the boundary torus that
bounds a compressing disk in \(H_i\) . We identify
\(T^2\) with the boundary of the standard solid torus in \(S^3\) in
such a way that \( \beta \) is identified with \(\gamma_1\) and 
\( \alpha \) has intersection number \(0\) with  \( \gamma_2 \).
 Connect \(x\) to
\(y\) by an embedded curve in
 \(T^2\) disjoint from \(\alpha\) and \(y\) to \(x\) by an embedded
 curve in \(T^2\) disjoint from \( \beta \). To obtain the knot, push
 the interior of the second curve into the solid torus, so that it
 becomes  disjoint from the second curve. We denote the resulting knot by
 \(K(p,q,r,s)\). Note that this identification is not unique ---
 different values of \(q,r\), and \(s\) may well produce the same knot.

The knot Floer homology of \(K(p,q,r,s)\) can be  computed
by the method of \cite{OS7}. Its total rank is always \(p\). To
find the Khovanov homology, we used the above method to produce a
planar diagram of \(K\). The diagram was then simplified (and, if
possible, identified) using {\it Knotscape}. Finally, the Khovanov homology was
computed using {\it KhoHo}. 

The table below contains a list of those \((1,1)\) knots which were
examined and   found to have property FK. The knots in the table are
of the form \(K(p,q,r,s)\) where \(11 \leq p \leq 17\). This is an
admittedly unscientific sample; it was chosen by sorting all the knots
with a given value of \(p\) by their Alexander polynomial, and then
selecting one representative for each Alexander polynomial. The idea
was to avoid duplicates, since the same knot (or its mirror image) is
usually represented by several different \(K(p,q,r,s)\)'s. The
selection was further narrowed by discarding knots with fewer than
\(11\) crossings, those which
appeared likely to be \( \delta\)-thin ({\it i.e.} \(\Delta_K(-1) = p
\)) and those whose Alexander polynomials suggested that they had
\(L\)-space surgeries. (These tended to be too big for {\it KhoHo} to
handle.) 

The entries in the table may be explained as follows. The first and
second columns identify the knot as \(K(p,q,r,s) \), and by its
{\it Knotscape} number. The next column shows the Alexander
polynomial. To save space, we have abbreviated by writing only the
coefficients of non-negative powers of
\(t\). For example, the first entry in the table indicates
an Alexander polynomial of \(-t^{-3} + 2 t^{-2} - 1 + 2 t^2 -t^3
\). The next column contains the \( \delta \)-polynomial, 
which is the Poincar{\'e} polynomial of
\(\hfk(K)\) with respect to the \(\delta\)-grading. For example
\(\hfk(K(11,3,3,3) \) has \(3\) generators with \(\delta\)-grading
\(-1\), and \(8\) with \(\delta\)-grading \(-2\).  The last column
shows \(\tau\), which in all cases is equal to \(s/2\). 

$$
\begin{array}{|l|c|l|l|c|}
\hline
\rule{0pt}{10pt} \text{Knot}  & \text{Knotscape \#} & \Delta_K(t) &  \delta \
\text{polynomial} & \tau  \\
\hline
\hline
K(11,3,3,2) & 11_{n} 19 & -1+0+2-1 & 3 \delta^{-1} + 8 \delta^{-2} &
-1  \\
\hline
K(13,5,1,2) & 11_n38 &+1+1-1 & 5 \delta^0 + 8 \delta ^{-1} & 0  \\
\hline
K(13,4,3,3) & 13_n192 & + 1-1+0+2-1 & 5 \delta^{-2} + 8 \delta ^{-3} & -2
\\ \hline
K(13,4,4,1) & 12_n725 & -3+2-1+0-1+1 & 9 \delta^{-4} + 4 \delta ^{-5}
& -5  \\ \hline
K(15,6,2,2) & 12_n121 & -1+0+1 & 8 \delta^0 + 7 \delta^{-1} & -1 
\\ \hline
K(15,3,4,2) & 12_n749 & -1 + 1 -1 + 1 & 11 \delta^{-1} + 4 \delta
^{-2} & -2 \\ \hline
K(15,5,4,2) & 12_n591 & -1+2-2+0+1 & 4 \delta ^4 + 11 \delta ^3 & 4  \\ \hline
K(15,3,8,1) &  15_n41127 & -3+2+1-2+1 & 7  \delta^{-1} + 8 \delta^{-2}
& -1  \\ \hline
K(15,4,5,2) & 12_n502 & -1 + 1 + 1-3 +2 & 3 \delta ^{-3} + 11 \delta
^{-4} & -4  \\ \hline
K(15,5,3,4) & 15_n4863 & -1 + 1-1+0+2 -1 & 7 \delta ^{-3} + 8
\delta^{-4} & -3  \\ \hline
K(17,7,1,2) &  11_n79 & -3+4-2 & 16 \delta^{1} + \delta^0 & 0  \\ \hline
K(17,6,2,2) & 15_n80764 & -1+1+1-1 & \delta^0 + 8 \delta^{-1} + 8
\delta^{-2} & 0  \\ \hline
K(17,7,2,1) & 11_n57 & +3-1-2+3-1 & 12 \delta^3 + 5 \delta ^2 & 3 
\\ \hline
K(17,3,4,2) & 14_n21882 & -1+0+2-2+1 & 13 \delta^0 + 4 \delta ^{-1} &
-1 \\
\hline
\end{array}
$$

\section{Knots Without Property FK}
\label{Sec:Examples2}
In this section, we give some examples of knots which are known not to
satisfy property FK. Initially, such examples were rather hard to come
by, but the appearance of {\it KhoHo} has made them substantially easier
 find. Although the size of the sample here is too small to support
 even an optimistic conjecture, it is worth noting that all the
 examples described here share the following properties.
\begin{itemize}
\item {\it Large bridge number:} Although this invariant is hard to
  compute, it appears that all examples have bridge number \( \geq
  4\).

\item {\it Large \(\khr\):} In all the examples, the rank of \( \khr
  \) is greater than the rank of \(\hfk\). 

\item {\it Torsion in \(\khr\):} All examples have \( \Z/2 \) torsion
  in their reduced Khovanov homology. For small knots, at least, this
  phenomenon is quite rare. \( \khr \) is free for all knots with
  fewer than 13 crossings, and there are only four 13 crossing knots which
  have \( \Z/2\) torsion in \( \khr\). (They are \(13_n3663\),
  \(13_n4587\), \(13_n4639\), and \(13_n5016\).) Two of
  these knots appear in the list below. It would certainly be
  interesting to determine if the other two have property FK. 
\end{itemize}
With these generalities taken care of, we turn to the
examples. 

\subsection{Torus knots}
The \((4,5)\) and \((4,7)\) torus knots do not have property FK. In
fact, it seems likely that this is the case for all \((p,q)\) torus
knots with \(p,q > 3\), but these are the only ones for which we
were able to determine \(\khr\). For these two knots, the knot Floer
homology is free with Poincar{\'e} polynomial
\begin{align*}
P_{\hfk} (T_{4,5}) & = t^6 + t^5u^{-1} + t^2u^{-2} + u^{-5} +
t^{-2}u^{-6} + t^{-5}u^{-11} + t^{-6}u^{-12} \\
P_{\hfk} (T_{4,7}) &  = 
t^9 + t^{8}u^{-1}  + t^{5}u^{-2} + t^{4}u^{-3} +
t^2u^{-4} -u^{-7} + t^{-2}u^{-8}  \\ & \qquad + t^{-4}u^{-11} + t^{-5}u^{12} +
t^{-8}u^{-17} + t^{-9}u^{-18}. 
\end{align*}
In particular, the ranks of \( \hfk \) are \(7\) and \(11\)
respectively.  In contrast, the ranks of \( \khr\) are \(9\) and
\(17\). Their Poincar{\'e} polynomials are
\begin{align*}
P_{\khr} (T_{4,5}) & = q^{12} + q^{16}u^2 + q^{18}u^3 + q^{18}u^4 +
q^{22}u^5 + q^{20}u^6 + q^{24}u^7 + q^{24}u^8 + q^{26}u^9 \\
P_{\khr} (T_{4,7}) & = q^{18} + q^{22}u^2 + q^{24}u^3 + q^{24}u^4 +
q^{28} u^{5} + q^{26}u^{6}+q^{30}u^7+2q^{30}u^{8}  \\ &
\qquad  + 2q^{32}u^{9}  +
q^{32}u^{10} + 2q^{36}u^{11} + q^{36}u^{12} + q^{38}u^{12} + q^{38}u^{13}
\end{align*}
In addition, \(\khr(T_{4,5}) \) has a \(\Z/2\) summand in degrees
\(q^{22}u^7\) and \(q^{28}u^{10}\), and \(\khr(T_{4,7}) \) has a \( \Z/2\)
summand in degrees \(q^{28}u^6\), \(q^{34} u^9\), \(q^{34}u^{10}\), and
\(q^{40}u^{13} \). (The author would like to thank Alexander
Shumakovitch for providing the information for the \((4,7) \) torus
knot.)

Torus knots are positive, so these examples automatically satisfy \(s
= 2 \tau \). 

\subsection{Cables of the trefoil}
The knot Floer homology of certain cabled knots was computed by Hedden
in \cite{Hedden}. Although it seems difficult to compute the Khovanov
homology of most cabled knots, the \((2,n)\) cables of the trefoil are 
small enough to be attacked directly using {\it KhoHo}. The
\((2,5)\) and \((2,7)\) cables of the positive trefoil each have
planar diagrams with \(13\) crossings. (Their {\it Knotscape} numbers
are \(13_n4639\) and \(13_n4587\).) Using the methods of
\cite{Hedden}, their knot Floer homologies can be seen to have ranks
\(5\) and \(7\), respectively. (The author would like to thank Matt
Hedden for sharing this fact.) More precisely, their Poincar{\'e}
polynomials are given by
\begin{align*}
P_{\hfk} (C_{2,5}T) & = t^4-t^3u^{-1}+u^{-2} - t^{-3} u^{-7} +
t^{-4}u^{-8} \\
P_{\hfk} (C_{2,7}T) & = t^5 - t^4u^{-1} + tu^{-2} -u^{-3} + t^{-1}u^{-4}
-t^{-4}u^{-9} + t^{-5}u^{-10} 
\end{align*}
On the other hand, their reduced Khovanov homologies both have rank
\(11\). Their Poincar{\'e} polynomials are 
\begin{multline*}
P_{\khr} (C_{2,5}T)  =q^8 + q^{12}u^2 + q^{14}u^3 + q^{14}u^4 +
     q^{18}u^5 + q^{18}u^6 \\ 
+ q^{20}u^7 + q^{20}u^8 + q^{22}u^9 + q^{26}t^{11} +
     q^{28}u^{12} 
\end{multline*}
and
\begin{multline*}
P_{\khr} (C_{2,7}T)  = q^{10} + q^{14}u^2 + q^{16}u^3 + q^{16}u^4 +
     q^{20}u^5 + q^{18}u^6 + \\ 
q^{22}u^7 + q^{22}u^8 + q^{24}u^9 + q^{28}u^{11} +
     q^{30}u^{12}
\end{multline*}
Both knots have \(\Z/2\) torsion in \(\khr\); for the \((2,5)\) cable
it is in degrees \(q^{16}u^6\), \(q^{18}u^7\), \(q^{22}u^9\), and \(
q^{24}u^{10}\), while for the \((2,7)\) cable it is in degrees 
\(q^{20}u^7\), \(q^{24}u^9\), and \( q^{26}u^{10} \). 
In both cases, it is easy to check that \(s = 2\tau\). 

More generally, it can be seen that for any odd \(n\geq 5\), the rank
of the knot Floer homology of the \((2,n)\) cable of the trefoil is
less than the rank of its \( \khr\). (The rank of \(\hfk\)  is
computed by \cite{Hedden}, while the the rank of \(\khr\) can be
bounded below using the skein exact sequence.) These knots have
positive diagrams, so they automatically satisfy \( s = 2 \tau \). 

\subsection{\((1,1)\) knots} 
He we include three other \((1,1)\) knots which do not have property
FK. In the notation of section~\ref{Subsec:(1,1)} these are
\(K(13,4,2,1)=16_n207543\), \(K(15,3,6,1)=16_n246032\), and
\(K(15,5,3,2)\), which {\it Knotscape} was unable to reduce to fewer than
\(17\) crossings. Below, we give individual details for each knot.

\noindent \(\mathbf{K(13,4,2,1)} \): 
\(\hfk\) is rank \(13\) with Poincar{\'e} polynomial  
\begin{multline*}
P_{\hfk}(K(13,4,2,1)) = t^4 u^{7} + 2t^3u^6 + t^2u^5 + tu^3 + tu^2 +
u^2  + t^{-1}u \\ + t^{-1}+ t^{-2}u+2t^{-3} + t^{-4}u^{-1}.
\end{multline*}
\noindent  
\( \khr\) is rank \(15\), with Poincar{\'e} polynomial
\begin{multline*}
P_{\khr}(K(13,4,2,1)) = q^6u^6 + q^4u^5 + q^2u^4 + 2 u^3 + u^2 +
q^{-2}u^2  \\ + q^{-2}u + q^{-4}u + q^{-2} + q^{-6} + q^{-6}u^{-1} +
q^{-6}u^{-2} + q^{-8}u^{-2} + q^{-8}u^{-3}. 
\end{multline*}
and \(\Z/2\) summands in degrees \(q^{-4}u^{0}, q^{-6}u^{-1},
q^{-10}u^{-3}\), and \(q^{-12}u^{-4}\). 
This knot has \(\tau = -1\) and \(s=-2\).

\noindent \(\mathbf{K(15,3,6,1)} \): 
\(\hfk\) is rank \(15\) with Poincar{\'e} polynomial  
\begin{equation*}
P_{\hfk}(K(15,3,6,1)) = t^2 u^3 + 2t u^3 + t u^2 + 4 u^2 + u + 2  +
t^{-1} + 2 t^{-1}u + t^{-2}u^{-1}. 
\end{equation*}
\( \khr\) is rank \(17\), with Poincar{\'e} polynomial
\begin{multline*}
P_{\khr}(K(15,3,6,1))  = 
 q^6u^2 + q^2 u + q^2 + 1 + 2u^{-1} + u^{-2} +
q^{-4}u^{-2} + q^{-2}u^{-3} + q^{-4}u^{-3}  \\  + q^{-4}u^{-4} +
q^{-6}u^{-4} + 2 q^{-6}u^{-5}  + q^{-8}u^{-6} + q^{-10}u^{-7} +
q^{-12}u^{-8}. 
\end{multline*}
and \(\Z/2\) summands in degrees \(q^4u^2\) and \(q^{-2}u^{-1}\).  
This knot has \(s = 2\tau  = 0\).

\noindent \(\mathbf{K(15,5,3,1)} \): 
\(\hfk\) is rank \(15\) with Poincar{\'e} polynomial  
\begin{multline*}
P_{\hfk}(K(15,5,3,1)) = t^4u^2 + 2 t^3u + t^2 + t^2 u + tu + 2  \\  +
u^{-1} + t^{-1}u^{-1}+ t^{-2}u^{-2} + t^{-2}u^{-4} +  2 t^{-3}u^{-5} +
t^{-4}u^{-6}
\end{multline*}
\noindent  
\( \khr\) is rank \(21\), with Poincar{\'e} polynomial
\begin{multline*}
P_{\hfk}(K(15,5,3,1)) = q^6 u^3 + q^4 u^2 + q^2 u + 3 + u^{-1} +
q^{-2}u^{-1} + q^{-2} u^{-2} \\ + q^{-4} u^{-2} + q^{-2} u^{-3}  +
q^{-6}u^{-3} + q^{-4} u^{-4} + q^{-6} u^{-4} \\ + q^{-6} u^{-5} +
q^{-8}u^{-5} + 2 q^{-8}u^{-6} + q^{-10}u^{-7} + q^{-12} u^{-8} +
q^{-14} u^{-9} 
\end{multline*}
and \(\Z/2\) summands in degrees \(q^4u^2, q^2u, q^{-2}u^{-1}\) and
\(q^{-4}t^{-2}\). This knot has \(s = 2\tau  =  0\).


\end{document}